\documentclass[a4paper,12pt]{article}
\usepackage{amsmath}
\usepackage{amssymb}
\usepackage[dvips]{graphicx}

\newtheorem{theorem}{Theorem}[section]

\newtheorem{definitionhead}[theorem]{Definition}

\begin{document}

\title{Remark to the paper Describing the set of words generated by
interval exchange transformation, posted 15 November 2007}
\author{A.Ya.~Belov \and A.L.~Chernyat'ev}
\date{}

\maketitle

\newpage
\pagestyle{plain}
\newpage

Let us call subdivision {\it good}, if 1) set corresponding to
each symbol is convex (i.e. interval or (semi)closed interval). 2)
If points $A$ and $B$ corresponds to the some color and interval
$(A,B)$ has discontinuity point, then $f(A)$ and $f(B)$ has
different color. Every subdivision can be further divided into
good subdivision, old superword can be obtained from new one by
gluing letters. Hence in the section ``Equivalence of the set of
uniformly recurrent words generated by piecewise-continuous
transformation to the set of words generated by interval exchange
transformation'' one can consider only good subdivision.

\end{document}